\documentclass{article}
\usepackage[english]{babel}
\usepackage{amsmath, amsthm, fancyhdr, hyperref}

\newcommand{\version}{ }

\usepackage{amsmath, amsthm, fancyhdr}
\usepackage{amscd}
\usepackage{amsfonts}
\usepackage{amssymb}
\usepackage{mathrsfs}
\usepackage{dutchcal}
\usepackage{url}
\usepackage{titlesec, eucal}
\titleformat{\section}[hang]{\large\normalfont\sc}{\thesection.}{18pt}{}{}
\titleformat{\subsection}[runin]{\normalfont\sc}{\thesubsection.}{3pt}{}
\usepackage{tocloft}

\renewcommand{\P}{\mathrm{P}}

\newcommand{\M}{\mathcal{M}}

\newcommand{\Z}{\mathbb{Z}}

\renewcommand{\O}{\mathcal{O}}

\DeclareMathOperator{\OG}{OG}
\DeclareMathOperator{\Mark}{Mark}

\renewcommand{\phi}{\varphi}
\renewcommand{\epsilon}{\varepsilon}

\def\eqref#1{(\ref{#1})}

\newcommand{\rk}{{\rm rk}}

\newcommand{\cntrct}                
{\hspace{2pt}\raisebox{1pt}{\text{$\lrcorner$}}\hspace{2pt}}

\numberwithin{equation}{section}

\newcommand{\arr}{\xrightarrow}

\theoremstyle{definition}

\newtheorem{proposition}[subsection]{\sc Proposition}
\newtheorem{theorem}[subsection]{\sc Theorem}

\theoremstyle{definition}
\newtheorem{definition}[subsection]{\sc Definition}

\newtheorem{conjecture}[subsection]{\sc Conjecture}

\newtheorem*{question}{\sc Question}

\theoremstyle{plain}

\def\endproof{{\small\textsquare}}
\def\proof{{{\sc Proof: }}}

\begin{document}

\begin{center}
{\Large\sc
Markman--Markushevich principle and Lagrangian fibrations on O'Grady six-manifolds\\[4mm]}
Rodion D\'eev

\end{center}

{\hfill \it to Prof. F. A. Bogomolov}

{\hfill \it on the occasion of his 80th birthday}

{\hfill \it with admiration}

\vspace{7mm}

{\small 
\hspace{0.1\linewidth}
\begin{minipage}[t]{0.7\linewidth}
{\sc Abstract:} To any compact holomorphically symplectic manifold and a nef BBF-isotropic class we associate a quasi-projective variety called the {\it Markman--Markushevich space} which parametrises the Lagrangian fibrations with this isotropic class modulo Tate--Shafarevich twists. We give a systematic interpretation of known results about Tate--Shafarevich twists of Lagrangian fibrations using this space, employ it to obtain a conjectural result for the last deformation type for which such a result is not yet proposed, and make some first steps towards its proof. As a byproduct, we generalise the results of Kumar on existence Shioda--Inose partners for abelian surfaces to a slightly more general class of K3 surfaces, and associate a canonical K3 surface with transcendental lattice isometric to $U^2\oplus\langle-2\rangle^2$ to a smooth cubic surface with a choice of a Schl\"afli double-six. Finally, we construct a noncompact six-manifold with an abelian fibration associated to a cubic surface with a Schl\"afli double-six, and give its interpretation in terms of the K3 surfaces we introduced.
\end{minipage}
}

\vspace{3mm}

\tableofcontents

\vspace{3mm}


\vspace{1mm}
\rule{75mm}{0.1pt}
\vspace{-3mm}

\newpage

\section{Introduction}

\subsection{} Holomorphically symplectic manifolds are a mysterious object. When Bogomolov (among others) discovered what is now called the Beauville--Bogomolov--Fujiki form on the second cohomology of an holomorphically symplectic manifold, it surprised him so much that he was led to a conclusion that higher-dimensional holomorphically symplectic manifolds do not exist at all, being mere products of K3 surfaces. As~stated, this is not true; but all of the known higher-dimensional examples (not much are known!) retain some familiarity to the lower-dimensional ones. The purpose of this note is to sum up a strain of extant research which concerns a particular kind of such familiarity, and to formulate what we call the ``Markman--Markushevich principle'' for Lagrangian fibrations. Informally speaking, it says that any Lagrangian fibration on a higher-dimensional holomorphically symplectic manifold may be decomposed locally over its base and then reglued in such a way that its total space becomes a moduli space of objects on a simpler thing. These ``things'' are known for the three of the four known deformation types: a K3 surface for $\mathrm{K}3^{[n]}$ type, abelian surface for $\mathrm{Kum}^n$ type, and cubic fourfold\footnote{In this case, the words ``moduli space of objects'' should be understood in a slightly heterodox fashion, although there is a way to make it less so \cite{_Bottini_}.} for $\OG_{10}$ type. We complement this by proposing a conjectural ``thing'' for the $\OG_6$ type, and present some evidence towards why this might be true.

\subsection{} The paper is deliberately written virtually without proofs. The reason for this is twofold. Firstly, the recent socio-technological development almost devalued the~proofs as such, and spending good time carving out them carefully seems unreasonable. Secondly, most likely, the author will be driven out of professional mathematics by a combination of economic and political pressures very soon, and at the same time is still convinced that the definitions and the general view explicated in this note are valuable enough to share with these who will be not. Perhaps these silly words would be made flesh by someone with better academic situation and citizenship.

\subsection{} The paper is organised as follows. In Section \ref{essential-section}, the reader would find a reminder about the theorems of Markman and Markushevich about Tate--Shafarevich twists of Lagrangian fibrations on $\mathrm{K}3^{[n]}$-type manifolds, and some geometric ideas behind them. The Definition \ref{essential-defn} of the essential Hodge structure and Markman--Markushevich space of a Lagrangian fibration allows to formulate analogous theorems for the other deformation types, even these hitherto unknown. Sections \ref{cubic-section} and \ref{lucky-section} develop these ideas in the case of six-dimensional O'Grady manifolds, culminating in the formulation of the main Conjecture \ref{main-conjecture}. By way of this, a couple of results about certain K3 surfaces of high Picard rank (namely, with transcendental lattice isometric to $U^2\oplus\langle-2\rangle^2$) are obtained. We call these surfaces the {\it lucky} ones (cf.~Proposition \ref{lucky-definition}). A new period map for the moduli space of cubic surfaces with a choice of Schl\"afli double-six, whose target is the period domain for the lucky K3 surfaces, is introduced (Definition \ref{cubic-Shioda-Inose} and Proposition \ref{defined-on-double-sixs}). In Section \ref{tensors-section}, we construct an abelian fibration similar to that on the noncritical part of an O'Grady six-manifolds, and conjecture (Conjecture \ref{lsv-og6-conjecture}) that it may be symplectically compactified, thereby yielding the desired construction. Finally, the Section \ref{cremona-section} concludes by linking this fibration to the family of K3 surfaces described above.

\section{\normalsize Essential Hodge structures and Markman--Mar\-ku\-she\-vich spaces}\label{essential-section}
\subsection{} The difficulty with which the holomorphically symplectic manifolds, even endowed with additional structure like a Lagrangian fibration, lend themselves to classification is in sharp contrast with the nature of the fibres of such fibrations, which, thanks to a beautiful lemma of Voisin, are abelian varieties even when the total space of the fibration is non-projective, and thus relatively well-understood. This, perhaps, was the motivation behind the following theorem of Markushevich:

\theorem[Markushevich, \cite{_Markushevich_}] A non-isotrivial Lagrangian fibration on a $\mathrm{K3}^{[2]}$-type manifold whose fibres are principally polarised is a Tate--Shafarevich twist of a Beauville--Mukai fibration of a genus 2 K3 surface. \endproof

\subsection{} The proof of Markushevich is especially transparent: any principally polarised abelian surface is either a Jacobian of genus 2 curve (and thus contains its copy, unique up to a translation), or a product of elliptic curves. This embedding may be performed in families, albeit up to an arbitrary vertical shift, so that one gets a universal genus 2 curve over a plane. The restriction of the symplectic form to this threefold defines an holomorphic characteristic foliation, and quotienting it out produces the germs of the desired K3 along some genus 2 curves. To glue these germs into a K3, one and indivisible, one has to get rid of the ambiguity with the vertical shifts, which amounts to saying ``up to Tate--Shafarevich twist’’. Unfortunately, this proof does not quite work for dimensions higher than four. The proof of Markman does.

\definition \label{essential-defn} Let $X \to B$ be an holomorphically symplectic manifold with a Lagrangian fibration. Let $\eta \in H^2(X, \Z)$ be the inverse image of $\O_B(1)$. The weight two $\Z$-Hodge structure $EH^2(X, \eta) = {\eta^\perp}{/\eta}$ is called the {\it essential Hodge structure} of the fibration. The quotient of the period domain of Hodge structures of this type by the subgroup of the monodromy group preserving $\eta$ is called the {\it Markman--Markushevich space} and denoted by $\Mark(X, \eta)$.

\subsection{} The essential Hodge structure is an invariant of the fibration modulo Tate--Sha\-fa\-re\-vich twist. Moreover, since its signature is $(2,b_2-4)$, the stabiliser of $\eta$ acts on the period domain properly discontinuously, hence $\Mark(X, \eta)$ is a quasi-projective variety. Its dimension equals $b_2(X)-4$, and it comes with a natural variation of weight two $\Z$-structures of K3 type---namely, $\eta^\perp/\eta$. SYZ conjecture for holomorphically symplectic manifolds which is now Engel--Mauri theorem \cite{_Engel--Mauri_} implies that the entire Markman--Markushevich space parametrises actual Lagrangian fibrations. Now the theorem of Markman (and Markushevich) may be stated this way.

\theorem[\cite{_Markman_}, Theorem 1.5] For a specific BBF-isotropic vector $\eta \in H^2(K3^{[g]}, \Z)$, the Beauville--Mukai construction yields a dominant generically finite map $bm : \mathcal{F}_g \to \Mark(\mathrm{K3}^{[g]}, \eta)$ from the moduli space of polarised genus $g$ K3 surfaces. Moreover, there is a pointwise Hodge isometry $h^\perp \arr{\sim} EH^2(bm(S,h),\eta)$ where $(S, h) \in \mathcal{F}_g$ is a polarised K3 surface, and $h^\perp \subset H^2(S, \Z)$ is the orthogonal to the polarisation. \endproof

\subsection{} The Beauville--Mukai construction is specific for $\mathrm{K3}^{[n]}$ type, but its cousins are known for other deformation types. For a linear system on an abelian surface, the {\it Debarre fibration} is a fibration into Prym varieties whose total space is an holomorphically symplectic variety of Kummer type. To a cubic fourfold, Laza, Sacc\`a, and Voisin associated a compactification of the relative intermediate Jacobian of its hyperplane section, which is an holomorphically symplectic manifold of O'Grady-10 type. Notice that the dimensions of the respective Markman--Markushevich spaces are $3$ and $20$, respectively, which coincides with the dimension of moduli spaces of abelian surfaces and cubic fourfolds, from which they are produced. This motivates the following two statements, which we label not as ``conjectures'', as they are folklore theorems among the experts (cf.~e.~g.~\cite{_Bottini_} regarding the second part of the statement).

\proposition \label{Markman-for-other-two-types} The Debarre construction yields a dominant generically finite map $deb \colon \mathcal{A}_2 \to \Mark(\mathrm{Kum}^n,\eta)$ from the moduli space of abelian surfaces of certain type of polarisation for a certain BBF-isotropic vector $\eta$. The Laza--Sacc\`a--Voisin construction yields the same with moduli space of cubic fourfolds and $\Mark(\OG_{10},\eta)$ for some $\eta \in H^2(\OG_{10},\Z)$. Both of them yield Hodge isometries $H^2(A,\Z) \supset h^\perp \cong EH^2(deb(A,h),\eta)$ and $H^4(V,\Z) \supset (h^2)^\perp \cong EH^2(lsv(V), \eta)$ (in the first case, $h$ stands for the polarisation on an abelian surface $A$, in the second, for the hyperplane section of a cubic fourfold $V$). \endproof

\subsection{} Perhaps it's no coincidence that the existence of the LSV compactification as a manifold of $\OG_{10}$-type was originally a conjecture of both Donagi--Markman and Markushevich--Tikhomirov. Thus it seems not entirely unreasonable to call the results of such type the ``Mar\-ku\-she\-vich--Markman principle''. We, instead of writing down its proof for the cases of Proposition \ref{Markman-for-other-two-types}, rather wonder:

\question What does Markman--Markushevich principle predict for $\OG_6$?

\vspace{1mm}
\rule{75mm}{0.1pt}
\vspace{-3mm}
\section{A natural candidate which somehow doesn't fit}\label{cubic-section}

\subsection{} Suppose we have a class of low-dimensional geometric objects with a moduli space $\M$ and an associated variation of weight one polarised $\Z$-Hodge structures $F \to \M$ with a ``construction'' which is a surjective, generically finite map $\M \arr{f} \Mark(\OG_6,\eta)$ for some isotropic $\eta \in H^2(\OG_6,\Z)$ that also induces a family of Hodge isometries $f_* \colon F_m \to EH^2(f(m),\eta)$.

\subsection{} We know that $b_2(\OG_6) = 8$, so $\dim\M=4$. The rank of the ``primitive middle cohomology'' of our conjectural object should equal $\rk~EH^2(\OG_g,\eta) = b_2-2 = 6$. The dimension of the base of the Lagrangian fibration in all cases is the projective dimension of the linear system on the object from which it stems. In our case, it is 3. To sum up, we look for something embedded (or otherwise mapped) into $\P^3$, with rank of primitive middle cohomology $6$, which has moduli space of dimension $4$... let's be honest. How can it be anything other than a cubic surface?

\subsection{} The problem is no less obvious: the Hodge structure of the middle cohomology of a cubic surface is all $H^{1,1}$, and cannot be meaningfully twisted into anything else, as the monodromy group $W(E_6)$ of the cubic surfaces acts irreducibly on $h^\perp \subset H^2$. Thus the cubic surface is only a part of the datum, and ought to be complemented by a kind of ``level structure'' (with zero-dimensional space of choices for a given surface, in order not to ruin the dimension of the moduli).

\subsection{} It is worthy to notice that a smooth cubic does have an associated Hodge structure with polarisation of type $(2,4)$. Namely, if the cubic is given by an equation $F(x_0, x_1, x_2, x_3) = 0$, it has an associated Hessian quartic $H = \det\left(\frac{\partial^2F}{\partial x_i\partial x_j}\right) = 0$. Its desingularisation is a K3 surface of Picard rank 16, and thus its transcendental lattice is a Hodge structure of weight two with signature $(2,4)$. To be more precise, the transcendental lattice is isometric to $U\oplus U(2) \oplus A_2(2)$, which has discriminant of order 48 \cite{_Dolgachev--Keum_}. In particular, it cannot be isometric to a $EH^2(\OG_6,\eta)$ for any nef isotropic $\eta \in H^2(\OG_6,\Z)$: the latter lattice is isometric to $U^3\oplus\langle-2\rangle^2$, and has discriminant of order 4. Moreover, \cite{_Mongardi--Rapagnetta_} prove that the divisibility of $\eta$ is one, hence the essential Hodge lattice is isometric to $U^2 \oplus \langle -2\rangle^2$. We shall call this lattice by $W$ for brevity.

\subsection{} The Hessian example, though unfit, is instructive. The lattice $W$ is a sublattice in the K3 lattice, and the moduli space of K3 surfaces with transcendental lattice isometric to $W$ yields a four-dimensional moduli space. The embedding of $W$ into the K3 lattice is unique thanks to the Nikulin's criterion.

\subsection{} A four-parametric family of K3 surfaces with a very similar transcendental lattice, namely, $U(2)^2 \oplus \langle-2\rangle^2$, was given by Matsumoto, Sasaki, and Yoshida \cite{_MSY_}. It is nothing more than the double cover of $\P^2$ ramified in a union of very general six lines. Its Picard rank is $16$, one class coming from $H^2(\P^2)$ and $15$ from the pairwise intersection of these lines. When the ramification locus is non-generic, the rank jumps (most famously, when the six lines are all tangent to a conic, in which case it yields a Kummer surface), so these are not of interest to us. To obtain a K3 surface with transcendental lattice $W$, one has to take an appropriate double cover of a Matsumoto--Sasaki--Yoshida K3 surface.

\vspace{1mm}
\rule{75mm}{0.1pt}
\vspace{-3mm}
\section{Shioda--Inose partners for cubic surfaces}\label{lucky-section}

\subsection{} Before embarking onto looking for the desired quotient, let us review what happens for generalised Kummer sixfolds through the same lens to ensure that we're not fooling ourselves. For a generalised Kummer manifold, $H^2$ with the BBF form is isometric to $U^3\oplus\langle-8\rangle$. They also admit Lagrangian fibrations with $(1,2,2)$-polarisation on a fibre, but unlike for the O'Grady sixfolds, the isotropic vector has divisibility $2$ \cite{_Floccari_,_Wieneck_}. Thus their essential $H^2$ is isometric to $U^2 \oplus \langle-2\rangle$ (isometric to $h^\perp \subset H^2(A, \Z)$ in line with Proposition \ref{Markman-for-other-two-types}).

\subsection{} This is the transcendental lattice of so-called Kumar K3s, which are double covers (or, dually, quotients) of Jacobian Kummer surfaces (they are also called Shioda--Inose, or sometimes Morrison, partners, of the corresponding abelian surfaces). Their other characterisation is having an elliptic fibration with two singular fibres of types $E_7$ and $E_8$ and a section (such a fibration, if exists, is unique up to an automorphism) \cite{_Kumar_}.

\subsection{} Let us describe the double covers of Kummer surfaces of Kumar type. First, let us construct another elliptic fibration. Name the rational curves on the Kumar surface so that the $E_7$ fibre becomes $e_1 + 2e_2 + 3e_3 + 4e_4 + 3e_5 + 2e_6 + e_7 + 2e_8$ ($e_4$ and $e_8$ are the classes that meet), and the $E_8$ fibre becomes $f_1 + 2f_2 + 3f_3 + 4f_4 + 5f_5 + 6f_6 + 4f_7 + 2f_8 + 3f_9$ ($f_6$ and $f_9$ meet). Let also $s$ be the section (it only intersects $e_1$ and $f_1$ out of the classes above). Opening the parentheses tells that the class $e_5 + e_8 + 2(e_4 + e_3 + e_2 + e_1 + s + f_1 + f_2 + f_3 + f_4 + f_5 + f_6) + f_7 + f_9$ has square zero, and it's a fibre of type $I^*_{10}$ in the elliptic fibration it defines. The curves $e_6$ and $f_8$ intersect it once, and hence are a pair of sections, while the class $e_7$ intersects the section $e_6$ but not $f_8$ (which implies that there exists another class $e_8$ such that $e_7 + e_8$ forms an $I_2$-type fibre, where $f_8$ intersects $e_8$). The other singular fibres are the six $I_1$s. The difference of the two sections is a 2-torsion element preserving the fibration, that is, a symplectic involtuion. It transforms the $I^*_{10}$ fibre into $I^*_5$ fibre, the $I_2$ fibre into an $I_1$ fibre, and six $I_1$ fibres into six $I_2$ fibres (the six of eight fixed points are on the $I_1$ fibres and two more are on the $I^*_{10}$ fibre). Its quotient is a Kummer surface.

\proposition \label{lucky-definition} The double plane model of a Kummer surface allows it to be deformed into a Matsumoto--Sasaki--Yoshida surface. Some fibration with a fibrewise double cover that yields the Kumar surface survives in this deformation. The covering surface has transcendental lattice isometric to $W = U^2 \oplus \langle -2\rangle^2$ as soon as its Picard rank is $16$. For all Kumar double covers that survive, all the covering surfaces are isomorphic. The double covers of Matsumoto--Sasaki--Yoshida K3 surfaces thus obtained are characterised by having an elliptic fibration with two $E_7$-type singular fibres and a section. We call such K3 surfaces the {\it lucky} ones.\\
\proof The double cover of a K3 is determined by an eight-tuple of disjoint rational curves in which it is ramified. For a double cover from the Kumar K3 to a Kummer K3, these eight curves are six components of the six $I_2$ fibres, and a pair of antennae of a $I^*_5$ fibre $E = h_0 + h_1 + 2(h_2 + h_3 + h_4 + h_5 + h_6 + h_7) + h_8 + h_9$. When a Kummer K3 undergoes a deformation into a generic Matsumoto--Sasaki--Yoshida K3, the classes $h_8, h_9$ of the two antennae cease to be algebraic, although its sum survives. In particular, $E$ may be now rewritten as the sum of algebraic classes $h_0 + h_1 + 2(h_2 + h_3 + h_4 + h_5 + h_6) + h_7 + (h_7 + h_8 + h_9)$. Let us proclaim the curves $h_7$ and $h_7 + h_8 + h_9$ the centres of ramification (the other six curves in which the ramification occurs are unaffected by the deformation). The corresponding double cover is a lucky K3 surface. Once its base degenerates back to a Kummer surface, the ramification divisor goes to $2h_7 + h_8 + h_9$ plus the six classes from $I_2$ fibres, and $2h_7$ may be safely subtracted once we take normalisation. Different choices of the fibration yield isomorphic Hodge structures, and hence isomorphic lucky K3 surfaces associated to six generic lines in a plane. \endproof

\subsection{} The elliptic fibration on a lucky K3 surface which admits an involution may be built the same way as that on a Kumar K3: two of the legs of the $E_n$-type singular fibres are joined through a section, which allows to isolate from this diagram a $D_m$-diagram (which serves as a fibre of the new fibration), and to use the rest of the $E_n$-diagrams as a pair of sections of that new fibration plus some component of fibres they meet elsewhere. The difference of the two sections is the desired symplectic involution.

\subsection{} A machine calculation indicates that after taking the involution, there is an appropriate choice of a class $H$ with $H^2=2$ and a blowdown of 15 (-2)-curves $E_{ij}$ so that the fibre of the new fibration becomes twice the class $5H - 3E_{12} - 2(E_{13} + E_{46} + E_{56}) - (E_{24} + E_{25} + E_{36} + E_{45})$. In other words, in the six-lines model the fibration becomes a pencil of rational quintics with a specific configuration of one triple point, three double points, and four ordinary points in the prescribed vertices of a complete hexagon. Then survival of the fibration under an arbitrary displacement of the six lines is evident, although the picture becomes painfully asymmetric in the six indices.

\definition \label{cubic-Shioda-Inose} Let $(C,L_1,\dots,L_6)$ be a cubic surface with a choice of six pairwise skew lines ({\it marked cubic} for brevity). Let $\{p_1,\dots,p_6\} \subset \P^2$ be the points that give $(C,L_1,\dots,L_6)$ as a blowup, and $S$ be the Matsumoto--Sasaki--Yoshida surface obtained as the double cover of the dual $\left(\P^2\right)^*$ in the union of six lines $p_1^*\cup\dots\cup p_6^*$. Its lucky cover is called the {\it Shioda--Inose partner} of the marked cubic $(C,L_1,\dots,L_6)$.

In particular, this gives a particular period map for marked cubics, which seems to be well-suited for the purposes of holomorphically symplectic geometry.

\vspace{3mm}

Now we are fully prepared to formulate our ultimate

\conjecture[Markman--Markushevich principle for $\OG_6$-manfiolds] \label{main-conjecture} There is a~dominant, generically injective map from the moduli space of marked cubics to the Markman--Markushevich space of O'Grady six-dimensional manifolds which establishes a Hodge isometry of the essential Hodge structure of the latter with the transcendental lattice of the Shioda--Inose partner of the former. When the cubic surface degenerates to a nodal and the six marked lines to six lines through the node, the~corresponding O'Grady six-manifold becomes the O'Grady's original example stemming from a Jacobian of genus 2 curve.

\subsection{} The moduli space of marked cubics carries a natural involution: for each five of six disjoint lines, there is a line which intersects these five but not the sixth, and these six lines form a complementary six-tuple. The six-tuples of points in $\P^2$ with this relation are called {\it Clebsch dual}. In the dual planes, two six-tuples of lines which are in Clebsch duality yield isomorphic double covers and hence isomorphic Shioda--Inose partners of the corresponding marked cubics:

\proposition \label{defined-on-double-sixs} Let $\{p_1,\dots,p_6\}, \{q_1,\dots,q_6\} \subset \P^2$ be two Clebsch dual configurations of points. Then the desingularisation of the double covers of the dual $\left(\P^2\right)^*$ ramified in $\cup_ip_i^*$ and $\cup_jq_j^*$ are isomorphic. Consequently, the period map for the marked cubics is invariant with respect to the involution.\\
\proof In fact, one ought not pass to the double covers in order to prove this. Consider a blowup of $\left(\P^2\right)^*$ in fifteen intersection points $e_{ij} = p_i^* \cap p_j^*$ of six lines (the~Matsumoto--Sasaki--Yoshida K3 surface is the double cover of this surface ramified in six disjoint $(-4)$-curves). The same surface admits another blowdown onto $\P^2$: one should replace the exceptional curves $E_{12}, E_{34}, E_{56}$ with the strict transform of the lines $e_{12}e_{34}, e_{34}e_{56}, e_{56}e_{12}$ (here $E_{ij}$ is the exceptional curve blown in lieu of $e_{ij}$). This is the standard Cremona transformation. One has to show that the lines $\{p_i^*\}$ project to a six-tuple of lines whose dual set of points is Clebsch dual to $\{p_i\}$. Choose the co-ordinates in which $E_{12} = (0 : 0 : 1), E_{34} = (0 : 1 : 0), E_{56} = (1 : 0 : 0)$. Then the lines $p_i^*$ have the equations $y/x = a$, $y/x = b$, $x = c$, $x = d$, $y = e$, $y = f$, respectively, and their Cremona transforms are $y/x = a^{-1}$, $y/x = b^{-1}$, $x = c^{-1}$, $x = d^{-1}$, $y = e^{-1}$, $y = f^{-1}$. For the points in the dual $\P^2$ this now amounts to the classical fact (due to Coble) that the Clebsch duality is a particular case of so-called Gale duality, cf. \cite[Example 5.12]{_Eisenbud--Popescu_}. \endproof





\vspace{1mm}
\rule{75mm}{0.1pt}
\vspace{-3mm}
\section{Abelian fibration associated to a marked cubic}\label{tensors-section}

In this section, we propose a candidate for the OG6 fibration. 

\vspace{-2mm}

\subsection{} The OG6 type manifolds all have a $(\Z/2)^8$-action which preserves the Lagrangian fibration whenever there is one \cite{_Mongardi--Rapagnetta_}. On its base $\P^3 \cong \P(V)$, this group acts through its quotient $(\Z/2)^4 = H$. This is the standard projective $G$-action which lifts to the four-dimensional linear representation of the Heisenberg group $H_{2,2}$ of order $64$ on $V$. This is the action of the 2-torsion subgroup $A[2]$ on the Kummer quartic embedding of $A/\{\pm 1\}$ \cite{_Eklund_}.

\vspace{-3mm}

\subsection{} \label{quadratic-tensor} The space of invariant quartic polynomials $W = \left(\mathrm{Sym}^4V^*\right)^{H_{2,2}}$ has dimension four. The Kummer quartics form an open dense subset of a divisor in $\P(W)$, namely, the Segre cubic threefold \cite{_Eklund_}. For a generic point $p \in \P(V)$, there is a unique $H$-invariant Kummer quartic with a node at $p$ (and hence with $Hp$ as its set of nodes). This defines an $H$-invariant meromorphic field $\mathcal{C} \subset \P T^*\P^3$ of quadratic cones on $\P(V)$: at each point, it is the tangent cone of the Kummer quartic which is singular at that point. This quadratic form is defined by the representation theory of $H_{2,2}$ alone.

\proposition Let $\tau \subset T_pP^3$ be a tangent plane at a generic point of a $\P^3$ with $H$-action. There exists a unique $H$-invariant Kummer surface containing $p$ as a smooth point with tangent plane $\tau$.\\
\proof The condition of passing through a point with precise tangent plane has codimension three, and thus defines a line in $\P(W)$ of all $H$-invariant quartics. It intersects the Kummer locus twice, once with multiplicity two in a point corresponding to the only invariant Kummer quartic singular at $p$, and has a residual intersection. \endproof

In particular, the Legendrian lifts of $H$-invariant Kummer quartics to $\P T^*\P^3$ form its rational Legendrian fibraion.

\vspace{3mm}

The classical construction of Cremona \cite{_Dolgachev_} associates to any marked cubic its representation as a section of the Segre cubic threefold. This motivates the following

\definition Let $\left(S, \{L_i\}_{i=1}^6\right)$ be a marked cubic. The pullback $\mathcal{Q}_S$ of the universal Kummer family along the Cremona section associated to this surface is called its {\it Cremona--Kummer family}. Its total space is a hypersurface in $\P T^*\P^3$ whose intersection with a generic fibre of the projection to $\P^3$ is a smooth plane cubic (a plane section of $S$).

\conjecture[LSV-type construction for an $\OG_6$-manifold] \label{lsv-og6-conjecture} Let $(S, \{L_i\})$ be a marked cubic, and $\mathcal{Q}_S \subset \P T^*\P^3$ its Cremona--Kummer family sitting inside the projectivised cotangent bundle of a $\P^3$ with $H$-action. Its preimage $\widetilde{\mathcal{Q}_S}$ inside the double cover of $\P T^*\P^3$ ramified in the conic bundle $\mathcal{C}$ defined in Paragraph \ref{quadratic-tensor} is a family of bielliptic curves of genus 4 through its projection to $\P^3$. The dual of the relative Prym of this family carries an holomorphic symplectic form, and admits a symplectic compactification of $\OG_6$ type.

\subsection{} \label{Matteini-paragraph} Another abelian fibration into Prym varieties of genus four double covers of elliptic curves is known: the Matteini fibration \cite{_Matteini_}. It takes a smooth cubic surface in $\P^3$ and takes its double cover ramified a canonical curve of genus $4$ which lies on it, and associates to its plane section the corresponding Prym variety. Its total space is an holomoprhic symplectic variety whose singularities have no symplectic resolution.

\section{Lagrangian origin of the Cremona--Kummer family}\label{cremona-section}

In this section, we give a description of the Cremona--Kummer families from the point of view of K3 surfaces discussed in the Sections \ref{cubic-section} and \ref{lucky-section}. Its purpose is to show that in the conjectural LSV construction of $\OG_6$-type manifolds the Cremona--Kummer family should play the r\^ole which in the original LSV construction is played by the cubic fourfold's Fano variety of lines.

\vspace{-2mm}

\subsection{} Let $Z \to \P^2$ be a double plane ramified in a sextic, so that $Z^{[2]}$ has a Lagrangian fibration over $\left(\P^2\right)^*$ whose fibre over the point $\ell \in \left(\P^2\right)^*$ is $\mathrm{Jac}\left(\tilde\ell \subset S\right)$. This works fine for generic branch sextic but once it splits off a line $\ell \subset \P^2$ the fibre gets bad above the dual point $\ell \in \left(\P^2\right)^*$, and we are going to have six of them. Let us blow up $\ell \in \left(\P^2\right)^*$: the points of the blown-into curve $L_\ell$ correspond to lines $a$ through $\ell$, that is, points $a \in \ell \subset \P^2$. Over the point $a \in C_\ell \subset \mathrm{Bl}_{\ell}\left(\P^2\right)^*$ let us place the Jacobian of the genus two curve formed as double cover of $\ell$ ramified in $a$ and the five other points of intersection of $\ell$ with the rest of the branching locus. This family of curves is not well-defined because all curves of genus two have an hyperelliptic involution, and the monodromy around the better part of the discriminant acts as it. Therefore, the extension of the Jacobian fibration to the blowup is impossible; but the Kummer fibration extends. Thus, any Matsumoto--Sasaki--Yoshida K3 gives rise to a Jacobian Kummer family $\mathcal{Q}'\to S$ parametrised by a cubic surface $(S,L_1, \dots, L_6) = \mathrm{Bl}_{\ell_1,\dots,\ell_6}\Pi^*$.

\subsection{} To relate this family to the Cremona--Kummer fibration, we wish to project it into $\P^3$ with an $H$-action in an equivariant way. Each individual member has an $H$-equivariant embedding, which is the map of the corresponding Jacobian with twice the theta-divisor. To glue them into one map from the total space of the entire family, the latter should be twisted.

\subsection{} \label{mukai-twist} The prototypical example of such a twist is the following. Let $Y \subset \P^5$ be a K3 surface with a polarisation of genus 5. Such a surface is always an intersection of three quadrics $Q_1, Q_2, Q_3$. To any number $m+1$ of quadrics in any projective space $\P^{n-1}$ one may associate their joint characteristic polynomial $\det\sum_{i=0}^m\lambda_iQ_i$. It defines a degree $n$ hypersurface in $\P^m$, which is in our case a sextic curve in the plane $\Pi$. Now consider a map $Y^{[2]} \to \Pi^*$: a length two subscheme of $S$ defines a line, and the quadrics in the net spanned by $Q_1, Q_2, Q_3$ which contain this line constitute a line in $\Pi$, that is, a point in $\Pi^*$. This fibration is a Tate--Shafarevich twist of the Jacobian fibration of the double cover of $\Pi$ ramified in the characteristic sextic.

\vspace{3mm}

We are indebted to Prof.~I.~Dolgachev for pointing out the relevant family of octic K3s associated to a marked cubic.

\definition A K3 surface is called {\it diagonal} if it may be represented as an intersection of three diagonal quadrics in $\P^5$.

\proposition The Hilbert square of the diagonal octic has an involution which acts by fibrewise $-1$ on the Lagrangian fibration, and whose quotient is birational to the Cremona--Kummer family of the corresponding marked cubic.\\
\proof For three diagonal quadrics, the joint characteristic polynomial splits into linear factors, so the partner genus 2 K3 of a diagonal octic $Y$ is a Matsumoto--Sasaki--Yoshida K3. The group generated by $x_i \mapsto \pm x_i$ preserves the diagonal octics, and acts as $(\Z/2)^5$ on them. Its subgroup which acts on the octic symplectaly is $(\Z/2)^4 = H$. Since the class of the Tate--Shafarevich twist between $\mathrm{Jac}\left(\tilde\ell \subset \widetilde{\Pi}\right)$ and $Y^{[2]} \to \Pi^*$ belongs to the 2-torsion of the Tate--Shafarevich group, it descends to the Kummer families. The quotient of $Y^{[2]}$ by a non-symplectic automoprhism is that Kummer fibration, and it may be birationally modified accordingly to create a family $\mathcal{Q}$ of Kummer surfaces over $S = \mathrm{Bl}_{L_1,\dots,L_6}\Pi^*$. The $H$-action preserves the fibres, and hence the $H$-equivariant maps from the Kummers to $\P^3$ glue into a single $H$-equivariant map $\mathcal{Q} \to \P^3$, so it is indeed its Cremona--Kummer family. \endproof

\proposition Take an elliptic fibration with two E7-type fibres on a lucky K3 surface. There exists a double cover of it ramified in the eight curves which form odd-multiplicity part of the two E7 fibres, and its blowdown is isomorphic to a diagonal octic K3 surface.\\
\proof By drawing a diagram, one may assure oneself that the E7 fibres of the lucky elliptic fibration turn into $I^*_0$-type fibres, whereas the six $I_1$ fibres turn into six $I_2$ fibres. The section of the lucky fibration turns into a bisection of this new fibration. One has to prove that a K3 surface with two $I^*_0$-type fibres, six $I_2$-type fibres, and a bisection is a diagonal octic. Notice that fibration with such singular fibres but with a section instead of a bisection is a Matsumoto--Sasaki--Yoshida K3 surface with a fibration given by the preimage of the pencil of lines from one of the nodes of the branching sextic (the two $I^*_0$ fibres lift from the lines passing through it, the six $I_2$ fibres lift from the six diagonals, and the section is the exceptional curve at the node). The twisted partnership which relates a genus 2 K3 surface with a section to a genus 5 K3 surface with a bisection is the same as in the Paragraph \ref{mukai-twist}, and the branching divisor being six lines yields a basis in which the three conics containing the genus 5 surface become diagonal. \endproof


\paragraph*{\normalfont\sc Acknowledgements.}

I am indebted to Anna Abasheva, Ignacio Barros, Simone Billi, Igor Dolgachev, Joel Fine, Grzegorz and Micha\l~Kapustka, Eyal Markman, Kieran O'Grady, \'Angel David R\'\i os Ortiz, Alexander Petrov, Andr\'es Rojas, Vasily Rogov, and Justin Sawon for interesting and fruitful discussions. Moreover, this note would have never been written without the hospitality and the research atmosphere of the University of Saarbr\"ucken.

\vspace{3mm}

\noindent {\sc Rodion D\'eev\\
Universit\'e libre de Bruxelles\\
CP 218, Boulevard du Triomphe,\\ 
Ixelles, Belgium} \\
\tt rodion.deev@ulb.be

\end{document}